\documentclass{article}
\newtheorem{theorem}{Theorem}[section]
\newtheorem{lemma}{Lemma}[section]

\newtheorem{definition}{Definition}[section]
\newtheorem{proposition}{Proposition}[section]

\newenvironment{proof}{{\bf Proof.}}{\par\hspace{25em}\rule{1ex}{1ex}\par}
\newcommand{\ds}{\displaystyle}
\title{On the intrinsic geometry of a unit vector field\footnote{Comment. Mat. Univ. Carolinae 43, 2 (2002), 299-317.}}
\author{Yampolsky A.}
\date{}
\begin{document}
\maketitle
\begin{abstract}
    We study the geometrical properties of a unit vector field on a Riemann\-ian 2-manifold,
    considering the field as a local
    imbedding of the manifold into its tangent sphere bundle with the Sasaki metric.
    For the case of constant curvature $K$, we give a description of the totally geodesic unit
    vector fields for $ K=0$ and $K=1$ and prove a non-existence result for $K\ne0,1$.
    We also found a family $\xi_\omega$ of  vector fields on the hyperbolic 2-plane $L^2$
    of curvature $-c^2$ which generate foliations on $T_1L^2$ with leaves
    of constant intrinsic curvature $-c^2$ and of constant extrinsic curvature~$-\frac{c^2}{4}$.   \\[1ex]
    {\it AMS subject class:} Primary 53B25, 53C42; Secondary 46E25, 20C20
\end{abstract}

\section*{Introduction}
      A unit vector field $\xi$ on a Riemannian manifold $M$ is called {\it holonomic}
    if $\xi$ is a field of normals of some family of regular
    hypersurfaces in $M$ and {\it non-holonomic} otherwise. The geometry of
    non-holonomic unit vector fields has been developed by A.Voss at the end of the 19-th
    century. The foundations of this theory can  be found in \cite{Am}. Recently, the
    geometry of a unit vector field has been considered from another point of view.
    Namely, let $T_1M$ be the  unit tangent sphere  bundle of $M$ endowed with the
    Sasaki metric \cite{S}. If $\xi$ is
    a unit vector field on $M$, then one may consider $\xi$ as a
    mapping $\xi : M \to T_1M $ so that the image $\xi (M) $  is a
    submanifold in $T_1M$ with the metric induced  from $T_1M$. So, one may apply the
    methods from the study of the geometry of submanifolds to determine geometrical characteristics of a
    unit vector field. For example, the unit vector field $\xi$ is said to be {\it  minimal}
    if $ \xi(M)$ is  of minimal volume with respect to the induced
    metric  \cite{G-Z}. A number of examples of locally minimal vector unit fields
    has been found (see \cite{ BX-V1, BX-V2, GD-V1}). On the other
    hand, using the geometry of submanifolds, we may find the Riemannian, Ricci or scalar
    curvature of a unit vector field using the second fundamental form of the
    submanifold $\xi(M)\in T_1M$ found in \cite{Ym}. In this paper
    we apply this approach to the  simplest case when the base space is
    2-dimensional and hence the submanifold $\xi(M)\in T_1M$ is a hypersurface.
    \vspace{1ex}

    {\bf Aknowlegement.} The author expresses his thanks to E.Boeckx for valuable
    remarks and the referee for careful reading and corrections.

    \section{The results}

    Let $\xi$ be a given unit vector field. Denote by $e_0$ a unit vector field such
    that $\nabla_{e_0}\xi=0$. Denote by $e_1$ a unit vector field, orthogonal to
    $e_0$, such that
        $$ \nabla_{e_1}\xi=\lambda\eta,$$
    where $\eta$ is a unit vector field, orthogonal to $\xi$. The function $\lambda$
    is a {\it signed} singular value of a linear operator $\nabla\xi : TM\to\xi^\perp$
    (acting as $(\nabla\xi)X=\nabla_X\xi$). Set
    $$
    \nabla_\xi\,\xi=k\,\eta,\quad \nabla_\eta\,\eta=\kappa\,\xi.
    $$
    The functions  $k$ and $\kappa$ are the {\it signed} geodesic curvatures of the integral
    curves of the fields $\xi$ and $\eta$  respectively. We prove that $\lambda^2=k^2+\kappa^2$.

    Denote the {\it signed} geodesic curvatures
    of the integral curves of the fields $e_0$ and $e_1$ as $\mu$ and $\sigma$
    respectively. Then
    $$
    \nabla_{e_0}\,e_0=\mu\,e_1,\quad \nabla_{e_1}\,e_1=\sigma\,e_0.
    $$
    The rientations of the frames $(\xi,\eta)$ and $(e_0,e_1)$ are independent.
    Set $s=1$ if the orientations are coherent and $s=0$ otherwise.

    The following result ( Lemma \ref{Forms} ) is a basic tool for the study.
    \vspace{1ex}

    {\it \
    Let $M$ be a 2-dimensional Riemannian manifold of Gaussian curvature $K$. The second
    fundamental form $\Omega$ of the submanifold $\xi(M)\subset T_1M$ is given by
    $$
   \Omega=
   \left[\begin{array}{cc}
           \ds -\mu\,\frac{\lambda}{\sqrt{1+\lambda^2}} & \ds {(-1)^{s+1}}\frac{K}{2} +\frac{e_0(\lambda)}{1+\lambda^2}\\[2ex]
           \ds {(-1)^{s+1}}\frac{K}{2} +\frac{e_0(\lambda)}{1+\lambda^2}& \ds
            e_1\left(\frac{\lambda}{\sqrt{1+\lambda^2}}\right)
   \end{array}\right].
   $$
   }
    \vspace{1ex}

    Using the formula for the sectional curvature of $T_1M^n$, we find an expression
    for the Gaussian curvature of $\xi(M^2)$ ( Lemma \ref{main}).
    \vspace{1ex}

    {\it
        The Gaussian curvature $K_\xi$ of a hypersurface $\xi(M)\in T_1M$ is given
        by
        $$
         \begin{array}{ll}
        \ds K_\xi=\frac{K^2}{4}+\frac{K(1-K)}{1+\lambda^2}&\ds +(-1)^{s+1}\frac{\lambda}{1+\lambda^2}
        e_0(K)+\\[2ex]
        &\ds\frac12
        \mu e_1\left(\frac{1}{1+\lambda^2}\right)-\left((-1)^{s+1}\frac{K}{2}+
        \frac{e_0(\lambda)}{1+\lambda^2} \right)^2,
        \end{array}
       $$
       where $K$ is the Gaussian curvature of $M$.
       }
       \vspace{1ex}

       As applications of these Lemmas, we prove the following theorems.
    \vspace{1ex}

    {\bf Theorem \ref{Tgfield}}{\it \
    Let $M^2$ be a Riemannian manifold of constant Gaussian curvature $K$.
    A unit vector field $\xi$ generating a totally geodesic submanifold in $T_1M^2$ exists if
    and only if  $K=0$ or $K=1$. Moreover,
    \begin{itemize}
        \item[(a)] if $K=0$, then $\xi$ is either a parallel vector field or moving along a
                    family of parallel geodesics with constant angle speed.
                    Geometrically, $\xi(M^2)$ is either $M^2$ imbedded isometrically
                    into $M^2\times S^1$ as a factor or a (helical) flat submanifold in
                    $M^2\times S^1$;
        \item[(b)] if $K=1$, then $\xi$ is a vector field on a standard sphere $S^2$ which is parallel
        along the meridians and  moving along the parallels  with a unit angle speed.
        Geometrically, $\xi(M^2)$ is a part of totally geodesic
        $RP^2$ locally isometric to sphere $S^2$ of radius 2 in $T_1S^2\stackrel{isom}{\approx}RP^3.$
    \end{itemize}
    }
    \vspace{2ex}

    {\bf Theorem \ref{Const}}{\it \
         Let $M^2$ be a space of constant Gaussian curvature $K$.
        Suppose that $\xi$ is a unit geodesic vector field on $M^2$. Then
        $\xi(M^2)$ has constant Gaussian curvature in one of the following cases:
        \begin{itemize}
            \item[(a)] $K=-c^2<0$
            and $\xi$ is a normal vector field for the family of horocycles on
              the hyperbolic 2-plane $L^2$ of curvature $-c^2$.
                In this case, $\ds K_\xi=-c^2$ and therefore $\xi(M^2)$
                is locally isometric the base space;
            \item[(b)] $K=0$ and $\xi$ is a parallel vector field on $M^2$. In this case
                $K_\xi=0$ and $\xi(M^2)$ is  also locally isometric to the base
                space;
            \item[(c)] $K=1$ and $\xi$ is any (local) geodesic vector field on the standard
                sphere $S^2$. In this case, $K_\xi=0$.
        \end{itemize}
        }

        {\bf Theorem \ref{Foli}}{\it \
    Let $L^2$ be a hyperbolic 2-plane of constant curvature $-c^2$. Then $T_1L^2$
    admits a hyperfoliation with leaves of constant intrinsic curvature $-c^2$ and of constant
    extrinsic curvature $-\frac{c^2}{4}$. The leaves are generated by unit vector fields
    making a constant angle with a pencil of parallel geodesics on $L^2.$
        }

    \section{Basic definitions and preliminary results}
    Let $(M,g)$ be an $(n+1)$ -- dimensional Riemannian manifold  with
    metric $g$. Let $\nabla$ denote the Levi-Civita connection on $M$.
    Then $\nabla_ X \xi$ is always orthogonal to $\xi$ and hence,
    $(\nabla\xi)X\stackrel{def}{=}\nabla_X\xi :T_pM \to \xi^\perp_p$ is a linear operator at each
    $p\in M$. We define an adjoint operator $ (\nabla\xi)^*\,X :\xi^\perp_p \to
    T_pM$  by
    $$
        \left< (\nabla\xi)^*X,Y\right>_g = \left< X,\nabla_Y\xi \right>_g.
    $$
    Then there is an orthonormal frame $e_0, e_1, \dots , e_n $ in $T_pM$ and an
    orthonormal frame $f_1, \dots , f_n $ in $\xi_p^\perp$ such that
    \begin{equation}\label{sing}
        (\nabla\xi)e_0=0 , \quad (\nabla\xi)e_\alpha =\lambda_\alpha f_\alpha,
        \quad (\nabla\xi)^*\,f_\alpha=\lambda_\alpha e_\alpha ,
        \qquad \alpha=1, \dots , n ,
    \end{equation}
    where $\lambda_1, \lambda_{2}, \dots  \lambda_n$ are real-valued functions.
    \begin{definition} The orthonormal frames satisfying (\ref{sing}) are called
    {\em singular frames} for the linear operator $(\nabla\xi)$ and the real valued functions
    $\lambda_1, \lambda_{2}, \dots \lambda_n $ are called the (signed) {\it singular
    values} of the operator $\nabla\xi$ with  respect to the singular frame.
    \end{definition}

    Remark that the sign of the singular value is definied up to the directions of the
    vectors of the singular frame.

    For each $\tilde X\in T_{(p,\xi)}TM$ there is a decomposition
    $$
        \tilde X= X_1^h+X_2^v
    $$
    where $(\cdot)^h$ and $(\cdot)^v$ are the horizontal
    and vertical lifts of vectors $X_1$ and $X_2$  from $T_pM$ to
    $T_{(p,\xi)}TM$.  The Sasaki metric is defined by the scalar product of the
    form
    $$
    \big<\big<\tilde X,\tilde Y \big>\big>=\big< X_1,Y_1 \big>+
                                             \big< X_2,Y_2\big>,
    $$
    where $\big<\cdot\,,\cdot\big>$ means the scalar product with respect to metric $g$.

    The following lemma has been proved in \cite{Ym}.
\begin{lemma}\label{L1}
    At each point $(p,\xi)\in \xi(M)\subset TM$ the vectors
    \begin{equation}\label{tang}
    \left\{
    \begin{array}{l}
        \displaystyle\tilde e_0  =  e_0^h, \\ [1ex]
        \displaystyle\tilde e_\alpha = \frac{1}{\sqrt{1+\lambda_\alpha^2}}(e_\alpha^h + \lambda_\alpha
        f_\alpha^v),
        \mbox{\hspace{3em}} \alpha=1,\dots , n ,
    \end{array} \right.
    \end{equation}
    form an orthonormal frame in the tangent space of $\xi(M)$ and the vectors
    \begin{equation}\label{norm}
    \tilde n_{\sigma |} =\frac{1}{\sqrt{1+\lambda_\sigma^2}}\big(-\lambda_\sigma
    e_\sigma^h +f_\sigma^v \ \big),\mbox{\hspace{3em}} \sigma=1,\dots , n ,
    \end{equation}
    form  an orthonormal frame in the normal space of $\xi(M)$.
\end{lemma}

    Let $R(X,Y)\xi=[\nabla_X,\nabla_Y]\,\xi-\nabla_{[X,Y]}\,\xi$ be the curvature tensor
    of $M$. Introduce the following notation
    \begin{equation}\label{r}
        r(X,Y)\xi=\nabla_X\nabla_Y\xi-\nabla_{\nabla_XY}\xi.
    \end{equation}
    Then, evidently,
    $$
    R(X,Y)\xi=r(X,Y)\xi-r(Y,X)\xi.
    $$

The following Lemma has also been proved in \cite{Ym}.

\begin{lemma}\label{L2}
    The components of second fundamental form of
    $\xi(M)\subset T_1M$ with respect to the frame (\ref{norm}) are given by
    $$
\begin{array}{rcl}
    \tilde \Omega_{\sigma | 00} &=&
    \frac{1}{\sqrt{1+\lambda_\sigma^2}}  \big< r(e_0,e_0)
    \xi,f_\sigma \big>, \\[3ex]
    \tilde \Omega_{\sigma | \alpha
    0} &=& \frac12 \frac{1}{\sqrt{(1+\lambda_\sigma^2)(1+\lambda_\alpha^2)}}
     \Big[ \big<
    r(e_\alpha,e_0) \xi + r(e_0,e_\alpha) \xi,f_\sigma \big> + \\[1ex]
    &&\hspace{3cm}  \lambda_\sigma \lambda_\alpha \big< R(e_\sigma,e_0) \xi, f_\alpha
    \big> \Big],   \\[2ex]
    \tilde \Omega_{\sigma | \alpha \beta}
    &=& \frac{1}{2}\frac{1}{\sqrt{(1+
    \lambda_\sigma^2)(1+\lambda_\alpha^2)(1+\lambda_\beta^2)}}
    \Big[ \big<
    r(e_\alpha, e_\beta) \xi+ r(e_\beta, e_\alpha) \xi, f_\sigma \big>\\[2ex]
    &&+ \lambda_\alpha \lambda_\sigma \big< R(e_\sigma, e_\beta)
    \xi, f_\alpha \big> +  \lambda_\beta \lambda_\sigma \big<
    R(e_\sigma, e_\alpha) \xi, f_\beta \big> \Big],
\end{array}
$$
    where $\{e_0,e_1,\dots,e_n;f_1,\dots,f_n\}$ is a singular frame of $(\nabla\xi)$
    and $\lambda_1,\dots,\lambda_n $ are the corresponding singular values.
\end{lemma}

    Let $\tilde\nabla$ and $\nabla$  be the Levi-Civita connections of the Sasaki metric of $TM$
    and the metric of $M$ respectively. The Kowalski formulas \cite{Kow} give the covariant
    derivatives of combinations of lifts of vector fields.
    \begin{lemma}[O.Kowalski]\label{Kow}
    Let $X$ and $Y$ be vector fields on $M$. Then  at each point $(p,\xi)\in TM$ we
    have
    $$
    \begin{array}{l}
    \tilde{\nabla}_{X^h}Y^h
    =(\nabla_XY)^h-\frac{1}{2}\left(R(X,Y)\xi\right)^v,  \\ [2ex]
    \tilde{\nabla}_{X^h}Y^v
    =\frac{1}{2}\left(R(\xi,Y)X\right)^h+(\nabla_XY)^v,  \\ [2ex]
    \tilde{\nabla}_{X^v}Y^h =\frac{1}{2}\left(R(\xi,X)Y\right)^h, \\[2ex]
    \tilde{\nabla}_{X^v}Y^v=0,
    \end{array}
    $$
    where $R$ is the Riemannian curvature tensor of $(M,g)$.
    \end{lemma}

    This basic result allows to find the curvature tensor of $TM$ (see \cite{Kow}) and
    the curvature tensor of $T_1M$ (see \cite{Bx-V3}). As a corollary, it is not
    too hard to find an expression for the {\it sectional curvature} of $T_1M$.
    It is well-known that $\xi^v$ is a unit normal for $T_1M$ as a hypersurface in $TM$.
    Thus, $\tilde X=X_1^h+X_2^v$ is tangent to $T_1M$ if and only if $\big<X_2,\xi\big>=0$.

    Let $\tilde X=X_1^h+X_2^v$ and $\tilde Y=Y_1^h+Y_2^v$, where $X_2,Y_2 \in \xi^\perp$,
    form an orthonormal base of a 2-plane $\tilde\pi \subset T_{(p,\xi)}T_1M$.
    Then we have \cite{Ym1}:
    \begin{equation}\label{Sec}
    \begin{array}{rl}
                \tilde K(\tilde\pi)=
     &\big< R(X_1,Y_1)Y_1,X_1\big>-\frac{3}{4}\|R(X_1,Y_1)\xi\|^2+\\ [2ex]
     & \frac{1}{4}\|R(\xi,Y_2)X_1+R(\xi,X_2)Y_1\|^2+\|X_2\|^2\|Y_2\|^2 -\big<X_2,Y_2\big>^2 + \\[2ex]
     & 3\big< R(X_1,Y_1)Y_2,X_2\big>-\big< R(\xi,X_2)X_1,R(\xi,Y_2)Y_1\big>+\\ [2ex]
     & \big< (\nabla_{X_1}R)(\xi,Y_2)Y_1,X_1\big> +\big< (\nabla_{Y_1}R)(\xi,X_2)X_1,Y_1\big>.
    \end{array}
    \end{equation}

    Combining the results of Lemma \ref{L1}, Lemma \ref{L2} and (\ref{Sec}), we can write an
    expression for the sectional curvature of $\xi(M)$.

    \begin{lemma}
      Let $\tilde X$ and $\tilde Y$ be an ortonormal vectors  which span a 2-plane $\tilde\pi$
    tangent to $\xi(M)\subset T_1M$. Denote by $K_\xi(\tilde\pi)$ the sectional curvature
    $\xi(M)$ with respect to metric, induced by Sasaki metric of $T_1M$. Then
    \begin{equation}\label{tensor}
    K_\xi(\tilde\pi)=\tilde K(\tilde\pi)+\sum_{\sigma}\left(\Omega_{\sigma|}(\tilde X,\tilde X)
    \Omega_{\sigma|}(\tilde Y,\tilde Y)-\Omega^2_{\sigma|}(\tilde X,\tilde Y)\right),
    \end{equation}
    where $\tilde K(\tilde\pi)$ is the sectional curvature
    of $T_1M$ given by (\ref{Sec}), $\Omega_{|\sigma}$ are the components of the second
    fundamental form of $\xi(M)$ given by Lemma \ref{L2} and the vectors  are
    given with respect to the frame (\ref{tang}).
    \end{lemma}

\section{The 2-dimensional case}

    Let $M$ be a 2-dimensional Riemannian manifold. The following proposition gives
    useful information about the relation between the singular values of the $(\nabla\xi)$-
    operator, geometric characteristics of the integral curves of singular frame and the
    Gaussian  curvature of the manifold.

    \begin{lemma}\label{prop1}
        Let $\xi$ be a given smooth unit vector field on $M^2$. Denote by $e_0$ a unit
        vector field on $M^2$ such that $\nabla_{e_0}\xi=0. $
        Let $\eta$ and $e_1$ be the unit vector fields on $M^2$ such that $(\xi,\eta)$
        and $(e_0,e_1)$ form two orthonormal frames on $M^2$. Denote by $\lambda$ a
        signed singular value of the operator $(\nabla\xi)$. Then we have
        $$ \nabla_{e_1}\xi=\lambda\eta,$$ and the following relations hold:

        \begin{itemize}
            \item[(a)] if $k=\big<\nabla_{\xi}\xi,\eta\big>$ is a signed geodesic curvature
        of a $\xi$-curve and $\kappa=\big<\nabla_{\eta}\eta,\xi\big>$ is a signed geodesic
            curvature of a $\eta$-curve, then

            $$\lambda^2=k^2+\kappa^2;$$

            \item[(b)] if $K$ is the Gaussian curvature of $M^2$, then
                        $$
                        (-1)^sK=e_0(\lambda)-\lambda\sigma,
                        $$
            where $\sigma=\big<\nabla_{e_1}e_1,e_0\big>$ is a signed geodesic curvature of a $e_1$-curve and
            $$s=\left\{\begin{array}{l} 1 \mbox{ if the frames $(\xi,\eta)$ and
            $(e_0,e_1)$ have the same orientation,} \\ 0 \mbox{ if the frames $(\xi,\eta)$ and
            $(e_0,e_1)$ have an opposite orientation}\end{array}\right..$$
        \end{itemize}
     \end{lemma}

    \begin{proof}
          (a) If $(\xi,\eta)$ is an orthonormal frame on $M^2$, then
           \begin{equation}\label{Fr1}
            \begin{array}{ll}
                \nabla_\xi\,\xi=k\,\eta,& \nabla_\xi\,\eta=-k\,\xi, \\[1ex]
                \nabla_\eta\,\xi=-\kappa\,\eta,& \nabla_\eta\,\eta=\kappa\,\xi.
            \end{array}
           \end{equation}
        Geometrically, the functions $k$ and $\kappa$ are the signed geodesic
        curvatures of $\xi$- and $\eta$-curves respectively.

        In a similar way we get
       \begin{equation}\label{Fr2}
            \begin{array}{ll}
                \nabla_{e_0}e_0=\mu\, e_1,& \nabla_{e_0}e_1=-\mu\, e_0, \\[1ex]
                \nabla_{e_1}e_0=-\sigma\, e_1,& \nabla_{e_1}e_1=\sigma\, e_0,
            \end{array}
        \end{equation}
           where $\mu$ and $\sigma$ are the signed geodesic curvatures of the
           $e_0$- and $e_1$-curves respectively.

        Let $\omega$ be an angle function between $\xi$ and $e_0$. Then we have two
        possible decompositions:
        $$
        \begin{array}{ll}
        \mbox{Or(+)  }
        \left\{\begin{array}{l}
                    e_0=\cos\omega \,\xi +\sin\omega\,\eta,\\
                    e_1=-\sin\omega\,\xi +\cos\omega\,\eta,
                \end{array}\right.
        & \mbox{Or(--)  }
        \left\{\begin{array}{l}
                    e_0=\cos\omega \,\xi +\sin\omega\,\eta,\\
                    e_1=\sin\omega\,\xi -\cos\omega\,\eta.
                \end{array}\right.
        \end{array}
        $$
        In the case $Or(+)$ we have
        $$
        \begin{array}{l}
        \nabla_{e_0}\,\xi=(k\cos\omega -\kappa\sin\omega)\,\eta ,\\
        \nabla_{e_1}\,\xi=-(k\sin\omega +\kappa\cos\omega)\,\eta,
        \end{array}
        $$
        and due to the choice of $e_0$ and $e_1$ we see that
        $$\left\{
        \begin{array}{lcl}
        k\cos\omega -\kappa\sin\omega&=&0, \\
        k\sin\omega +\kappa\cos\omega &=&-\lambda.
        \end{array}
        \right.
        $$
        So, for the case of $Or(+)$ \ $ k=-\lambda \sin\omega,\ \kappa =-\lambda\cos\omega.$

         In a similar way, for the case of $Or(-)$ \
         $k=\lambda \sin\omega,\  \kappa =\lambda\cos\omega.$
        In both cases
        $$\lambda^2=k^2+\kappa^2.$$

        (b) Due to the choice of the frames,
        $$
        \begin{array}{rl}
        \big<R(e_0,e_1)\xi,\eta\big>=&\big<\nabla_{e_0}\nabla_{e_1}\xi-\nabla_{e_1}\nabla_{e_0}\xi-
        \nabla_{\nabla_{e_0}e_1-\nabla_{e_1}e_0}\xi,\eta\big>=
        \\[1ex]
        &\big<\nabla_{e_0}(\lambda\,\eta)-\nabla_{-\mu\,e_0+\sigma\,e_1}\,\xi,\eta\big>=e_0(\lambda)-
        \lambda\sigma.
        \end{array}
        $$

        On the other hand,
       \begin{equation}\label{R}
        \big<R(e_0,e_1)\xi,\eta\big>=\left\{
        \begin{array}{l}
        -K \mbox{ for the case of $Or(+)$}, \\
        +K \mbox{ for the case of $Or(-)$.}
        \end{array}\right.
       \end{equation}
        Set $s=1$ for the case $Or(+)$ and $s=0$ for the case $Or(-)$.
        Combining the results, we get $(-1)^sK=e_0(\lambda)-\lambda\sigma,$
        which completes the proof.

    \end{proof}

   The result of Lemma \ref{L2} can also be simplified in the following way.

    \begin{lemma}\label{Forms}
    Let $M$ be a 2-dimensional Riemannian manifold of Gaussian curvature $K$. In terms of
    Lemma \ref{prop1}  the second
    fundamental form of the submanifold $\xi(M)\subset T_1M$ can be presented in two
    equivalent forms:
    \begin{itemize}
    \item[(i)]
    $$
   \Omega=
   \left[\begin{array}{cc}
           \ds -\mu\,\frac{\lambda}{\sqrt{1+\lambda^2}} & \ds {(-1)^{s+1}}\frac{K}{2} +\frac{e_0(\lambda)}{1+\lambda^2}\\[2ex]
           \ds {(-1)^{s+1}}\frac{K}{2} +\frac{e_0(\lambda)}{1+\lambda^2}& \ds
            e_1\left(\frac{\lambda}{\sqrt{1+\lambda^2}}\right)
   \end{array}\right],
   $$
       \item[(ii)]
    $$
   \Omega=
   \left[\begin{array}{cc}
           \ds -\mu\,\frac{\lambda}{\sqrt{1+\lambda^2}} & \ds \frac12\left(\sigma\,\lambda +\frac{1-\lambda^2}{1+\lambda^2}e_0(\lambda)\right)\\[2ex]
           \ds \frac12\left(\sigma\,\lambda +\frac{1-\lambda^2}{1+\lambda^2}e_0(\lambda)\right)& \ds e_1\left(\frac{\lambda}{\sqrt{1+\lambda^2}}\right)
   \end{array}\right].
   $$
    \end{itemize}

   \end{lemma}

   \begin{proof}

    At each point $(p,\xi)\in \xi(M)$ the vectors
    $$
    \left\{
    \begin{array}{l}
    \displaystyle\tilde e_0  =  e_0^h, \\ [1ex]
    \displaystyle\tilde e_1 = \frac{1}{\sqrt{1+\lambda^2}}(e_1^h + \lambda
        \eta^v)
    \end{array} \right.
    $$
    form an orthonormal frame in the tangent space of $\xi(M)$ and
    $$
  \quad  \tilde n =\frac{1}{\sqrt{1+\lambda^2}}\big(-\lambda e_1^h +\eta^v \ \big),
    $$
    is a unit normal for $\xi(M)\subset T_1M$.

    Thus we see that in a 2-dimensional case the components of $\Omega$ take the form
   $$
   \begin{array}{l}
    \ds \Omega_{00}=\frac{1}{\sqrt{1+\lambda^2}}\big<r(e_0,e_0)\xi,\eta\big>,
    \quad \ds\Omega_{11}=\frac{1}{(1+\lambda^2)^{3/2}}\big<r(e_1,e_1)\xi,\eta\big>,
    \\[3ex]
    \ds\Omega_{01}=\frac{1}{2}\frac{1}{1+\lambda^2}\Big[\big<r(e_1,e_0)\xi+r(e_0,e_1)\xi,\eta\big>+\lambda^2
    \big<R(e_1,e_0)\xi,\eta\big>\Big].
    \end{array}
   $$
   Keeping in mind (\ref{r}), (\ref{Fr2}) and (\ref{R}), we see that
    $$
    \begin{array}{ll}
    \ds\big<r(e_0,e_0)\xi,\eta\big>=-\mu\lambda, & \ds
    \big<r(e_0,e_1)\xi,\eta\big>=e_0(\lambda),\\[2ex]
    \ds\big<r(e_1,e_0)\xi,\eta\big>=\sigma\lambda,&
    \ds \big<r(e_1,e_1)\xi,\eta\big>=e_1(\lambda),\\[2ex]
    \big<R(e_0,e_1)\xi,\eta\big>=(-1)^sK.&
    \end{array}
    $$
  So we have
    $$
    \begin{array}{l}
    \ds \Omega_{00}= -\mu\,\frac{\lambda}{\sqrt{1+\lambda^2}},\qquad
   \ds \Omega_{11}=\frac{e_1(\lambda)}{(1+\lambda^2)^{3/2}}=e_1\left(\frac{\lambda}{\sqrt{1+\lambda^2}}\right),\\[3ex]
   \ds \Omega_{01}=\frac{1}{2(1+\lambda^2)}(e_0(\lambda)+\lambda\sigma-\lambda^2(-1)^sK)=
        \left\{\begin{array}{l}
                        \ds {(-1)^{s+1}}\frac{K}{2}
                        +\frac{e_0(\lambda)}{1+\lambda^2}\\[2ex]
                        \ds \frac12\left(\sigma\,\lambda +\frac{1-\lambda^2}{1+\lambda^2}e_0(\lambda)\right)
                \end{array}
        \right. ,
    \end{array}
    $$
    where Lemma \ref{prop1} (b) has been applied in two ways.
    \end{proof}

\subsection{Totally geodesic vector fields}
    The main goal of this section is to prove the following theorem.

    \begin{theorem}\label{Tgfield}
    Let $M^2$ be a Riemannian manifold of constant Gaussian curvature $K$.
    A unit vector field $\xi$ generating a totally geodesic submanifold in $T_1M^2$ exists if
    and only if $K=0$ or $K=1$. Moreover,
    \begin{itemize}
        \item[a)] if $K=0$, then $\xi$ is either a parallel vector field or is moving along a
                    family of parallel geodesics with constant angle speed.
                    Geometrically, $\xi(M^2)$ is either $M^2$ imbedded isometrically
                    into $M^2\times S^1$ as a factor or a (helical) flat submanifold in $M^2\times
                    S^1$;
        \item[b)] if $K=1$, then $\xi$ is a vector field on a sphere $S^2$ which is parallel
        along the meridians and  moving along the parallels  with a unit angle speed.
        Geometrically, $\xi(M^2)$ is a part of totally geodesic $RP^2$ locally isometric to
        sphere $S^2$ of radius 2 in $T_1S^2\stackrel{isom}{\approx}RP^3$ .
    \end{itemize}
    \end{theorem}

    The proof will be divided into a series of separate propositions.

    \begin{proposition}\label{Semigeo}
    Let $M^2$ be a Riemannian manifold. Let $D$ be a domain in $M^2$ endowed with
    a semi-geodesic coordinate system such that  $ds^2=du^2+f^2\, dv^2,$
    where $f(u,v)$ is some non-vanishing function.  Denote by
    $(e_0, e_1)$ an orthonormal frame in $D$ and specify $e_0=\partial_u\,, \,e_1=f^{-1}\partial_v$.
    If $\xi$ is a unit vector field in $D$ parallel along $u$-geodesics, then $\xi$ can be
    written given as  $$\xi=\cos \omega \,e_0+\sin \omega\,e_1,$$
    where $\omega=\omega(v)$ is an angle function and

    (a) a singular frame for $\xi$  may be chosen as
    $\big\{ e_0,\ e_1, \eta=-\sin \omega\,e_0+\cos\omega \,e_1\big\};$

    (b) a singular value for $\xi$ in this case is $\lambda =e_1(\omega)-\sigma,$
    where $\sigma$ is a signed geodesic curvature of the $e_1$-curves.

     \end{proposition}

    \begin{proof}
    Indeed, if $\xi$ is parallel along $u$-geodesics, then evidently the angle function
    $\omega$ between $\xi$ and the $u$-curves does not depend on $u$. So this function has the form
    $\omega=\omega(v)$ and $\xi=\cos \omega \,e_0+\sin \omega\,e_1$.
     Moreover, since
    $$
    \begin{array}{ll}
    \nabla_{e_0}e_0=0, & \nabla_{e_0}e_1=0, \\[1ex]
    \nabla_{e_1}e_0=\frac{f_u}{f}e_1, & \nabla_{e_1}e_1=-\frac{f_u}{f}e_0,
    \end{array}
    $$
    we see that $\sigma =-\frac{f_u}{f}$ and $\nabla_{e_1}\xi=(e_1(\omega)-\sigma)\,\eta,$
    where $\eta=-\sin\omega\,e_0+\cos\omega\,e_1$.
    Therefore, $\lambda=e_1(\omega)-\sigma$ and the proof is complete.
    \end{proof}

    \begin{proposition}\label{Neg}
    Let $M^2$ be a Riemannian manifold of constant negative curvature $K=-r^{-2}<0$. Then
    there is no totally geodesic unit vector field on $M^2$.
    \end{proposition}

    \begin{proof}
    Suppose $\xi$ is totally geodesic unit vector field on $M^2$.
    Set $\Omega\equiv 0$ in Lemma \ref{Forms}. Then $\lambda\mu\equiv 0$. If
    $\lambda\equiv 0$ in some domain $D\subset M^2$, then $\xi $ is parallel in this
    domain and hence $M^2$ is flat in $D$, which contradicts the hypothesis.
    Suppose that $\mu\equiv 0$ at least in some domain $D\subset M^2$.
    This means that $e_0$-curves are geodesics in $D$ and the field $\xi$ is parallel
    along them. Choose a family of $e_0$-curves and the orthogonal trajectories as a local
    coordinate net in $D$.
    Then the first fundamental form of $M^2$ takes the form
    $$
    ds^2=du^2+f^2\,dv^2,
    $$
    where $f(u,v)$ is some function. Since $M^2$ is of constant curvature $\ds K=-\frac{1}{r^2}$,
    the function $f$ satisfies the equation
    $$
    f_{uu}-\frac{1}{r^2}f=0.
    $$
    The general solution of this equation is
    $$
    f(u,v)=A(v)\cosh(u/r)+B(v)\sinh(u/r).
    $$
    There are two possible cases:
    $$
    \begin{array}{rl}
        (i)& A^2(v)\equiv B^2(v) \mbox{ over the whole domain $D$};\\[1ex]
        (ii)& A^2(v)\ne B^2(v)   \mbox{ in some subdomain $D'\subset D$}.
    \end{array}
    $$
    Case (i). In this case, in dependence of the signs of $A(v)$ and $B(v),$
    $$
    f(u,v)=A(v)e^{u/r}\quad \mbox{or} \quad f(u,v)=A(v)e^{-u/r}.
    $$
    Consider the first case (the second case can be reduced to the first one after the parameter
    change $u\to -u$). Making an evident $v$-parameter change, we  reduce the metric to the form
    $$
    ds^2=du^2+r^2e^{2\,u/r}\,dv^2.
    $$
    Applying Proposition \ref{Semigeo} for $f=re^{u/r}$, we get
    $\ds
     \lambda=\frac{1}{r}(\omega'\, e^{-u/r}+1).
    $
    Setting $\Omega_{11}\equiv 0$, we see that $e_1(\lambda)\equiv 0$. Hence $ \omega''=0$,
    i.e., $\omega=av+b$.
    Therefore,
    $$
    \lambda=\frac1r \left( a\, e^{-u/r}+1\right).
    $$
    Considering $\Omega_{01}\equiv 0$ (with $s=1$ because of $Or(+)$-case), we get
    $$
    -\frac{1}{2r^2}+\frac{\frac1re_0( a\, e^{-u/r}+1)}{1+\frac{1}{r^2}( e^{-u/r}a+1)^2}=
    -\frac{(\frac{1}{r^2}+1)(ae^{-u/r}+1)^2-a^2 e^{-2u/r}}{2\,r^2[ 1+\frac{1}{r^2}(ae^{-u/r}+1)^2]}\not\equiv
    0,
    $$
    and hence, this case is not possible.
    \vspace{1ex}

    Case (ii). Choose a subdomain $D'\subset D$ such that $A^2(v)<B^2(v)$ or $A^2(v)>B^2(v)$
    over $D'$. Then    the function $f$ may be presented respectively  in two forms:
    $$
    \begin{array}{rl}
     \mbox{(a)}& f(u,v)=\sqrt{B^2-A^2}\,\sinh(u/r+\theta) \mbox{\quad or
     }\\[1ex]
     \mbox{(b)}& f(u,v)=\sqrt{A^2-B^2}\cosh(u/r+\theta),
    \end{array}
    $$
    where $\theta(v)$ is some function.

    Consider the case (a).   After a $v$-parameter change, the
    metric in $D'$ takes the form
    $$
    ds^2=du^2+r^2\sinh^2(u/r+\theta)\,dv^2.
    $$
    Applying Proposition \ref{Semigeo} for $f=r\sinh(u/r+\theta)$, we get
    $$
      \lambda=\frac{\omega'}{r\,\sinh(u/r+\theta)} +\frac1r\coth(u/r+\theta).
    $$
    Considering $\Omega_{11}\equiv 0$, we have $e_1(\lambda)\equiv 0$ which implies
    the identity
    $$
    \omega''\,\sinh(u/r+\theta)-\omega'\theta'\,\cosh(u/r+\theta)-\theta'\equiv 0.
    $$
    From this we get $ \omega''=0,\ \ \theta'=0$ and hence
    $\ds
    \left\{
    \begin{array}{l}
        \theta=const, \\
        \omega=av+b
    \end{array}
    \right.
    $   ($a,b=const$).
    After a parameter change we reduce the metric to the form
    $$
    ds^2=du^2+ r^2\sinh^2(u/r)\,dv^2
    $$
    Applying Proposition \ref{Semigeo} for $f=r\sinh(u/r)$, we get
    $\ds
    \lambda = \frac{a+\cosh(u/r)}{r\sinh(u/r)}.
    $
    The substitution into $\Omega_{01}$ gives
    $$
    -\frac12\frac{(\frac{1}{r^2}+1)[a+\cosh (u/r)]^2-a^2+1}{r^2\sinh^2(u/r)+[a+\cosh
    (u/r)]^2}
    \not\equiv 0,
    $$
    which completes the proof for the  polar case.

    The {\it Cartesian} case  consideration gives
    $\ds \omega=av+b, \quad\lambda = \frac{a+\sinh(u/r)}{r\cosh(u/r)}$
    and
    $\ds
    \Omega_{01}=-\frac12\frac{(\frac{1}{r^2}+1)[a+\sinh (u/r)]^2-a^2-1}{r^2\cosh^2(u/r)+[a+\sinh
    (u/r)]^2}
    \not\equiv 0,
    $
    which completes the proof.
    \end{proof}

    \begin{proposition}\label{Pos}
    Let $M^2$ be a Riemannian manifold of constant positive curvature $K=r^{-2}>0$. Then
    a totally geodesic unit vector field $\xi$ on $M^2$ exists if $r=1$ and $\xi$ is
    parallel along the meridians of $M^2$ locally isometric to $S^2$ and moves along
    the parallels with a unit angle speed. Geometrically, $\xi(M^2)$ is a part of totally
    geodesic $RP^2$ locally isometric to sphere $S^2$ of radius 2 in
    $ T_1S^2\stackrel{isom}{\approx}RP^3$ .
    \end{proposition}

    \begin{proof}
    Suppose $\xi$ is totally geodesic unit vector field on $M^2$.
    The same arguments as in Proposition \ref{Neg} lead to the case $\mu\equiv0$ at
    least in some domain $D\subset M^2$. So, choose again a family
    of $e_0$-curves and the orthogonal trajectories as a local coordinate net in $D$.
    Then the first fundamental form of $M^2$ can be expressed as
    $\ds
    ds^2=du^2+f^2\,dv^2,
    $
    where $f(u,v)$ is some function. Since $M^2$ is of constant curvature $ K=r^{-2}$,
    the function $f$ satisfies the equation
    $$
    f_{uu}+\frac{1}{r^2}f=0.
    $$
    The general solution of this equation
    $\ds
    f(u,v)=A(v)\cos(u/r)+B(v)\sin(u/r)
    $
    may be presented in two forms:
     $$
    \begin{array}{rl}
     \mbox{(a)}& f(u,v)=\sqrt{A^2+B^2}\,\sin(u/r+\theta) \mbox{\quad or }\\
     \mbox{(b)}& f(u,v)=\sqrt{A^2+B^2}\cos(u/r+\theta),
    \end{array}
    $$
    where $\theta(v)$ is some function.

    Consider first, the case (a).
    After $v$-parameter change, the metric in $D$ takes the form
    $$
    ds^2=du^2+r^2\sin^2(u/r+\theta)\,dv^2.
    $$
    Applying Proposition \ref{Semigeo} for $f=r\sin (u/r+\theta)$, we get
    $$
     \lambda=\frac{\omega'}{r\,\sin(u/r+\theta)}+\frac{1}{r}\cot(u/r+\theta).
    $$
    Setting $\Omega_{11}\equiv 0$, we find $e_1(\lambda)\equiv 0$ which implies
    the identity
    $$
    \omega''\,\sin(u/r+\theta)-\omega'\theta'\,\cos(u/r+\theta)+\theta'\equiv 0.
    $$
    From this $ \omega''=0,\ \ \theta'=0$ and we have again
    $
    \left\{
    \begin{array}{l}
        \theta=const, \\
        \omega=av+b
    \end{array}
    \right.
    $
    \ $ a,b=const$.
    After a suitable $u$-parameter change, we reduce the metric to the form
    $$
    ds^2=du^2+ r^2\sin^2(u/r)\,dv^2
    $$
    Applying Proposition \ref{Semigeo} for $f=r\sin(u/r)$, we get
    $\ds
    \lambda = \frac{a+\cos(u/r)}{r\,\sin(u/r)}.
    $
    Substitution into $\Omega_{01}$ gives
    $$
    \frac12\frac{(\frac{1}{r^2}-1)[a+\cos (u/r)]^2+a^2-1}{r^2\sin^2 (u/r)+[a+\cos (u/r)]^2}\equiv 0,
    $$
    which is possible only if $r=1$ and $|a|=1$.
    So, we obtain to the standard sphere metric
    $$
    ds^2=du^2+\sin^2u \,dv^2
    $$
    and ( after the $\pm v+b \to v $ parameter change ) the unit vector field
    $$
    \xi=\left\{\cos v,\frac{\sin v}{\sin u}\right\}.
    $$
    This vector field is parallel along the meridians of $S^2$ and
    moves helically along the parallels of $S^2$ with unit angle speed.

     For the case (b) one can find
    $\ds
    \omega=av+b,\quad    \lambda=\frac{a-\sin(u/r)}{r\,\cos(u/r)}
    $
    and
    $$
    \Omega_{01}=\frac12\frac{(\frac{1}{r^2}-1)[a-\sin (u/r)]^2+a^2-1}{r^2\cos^2 (u/r)+[a-\sin (u/r)]^2}\equiv
    0,
    $$
    which gives $r=1$ and $|a|=1$ as a result. Thus, we have a  metric
    $$
    ds^2=du^2+\cos^2 u \,dv^2
    $$
    and a vector field
    $\ds
    \xi=\left\{\cos v,\frac{\sin v}{\cos u}\right\}.
    $
    It is easy to see that the results  of cases (a) and (b) are geometrically
    equivalent.

    Introduce the local coordinates $(u,v,\omega) $ on $T_1S^2$, where $\omega$ is the
    angle between arbitrary unit vector $\xi$ and the coordinate vector field
    $X_1=\big\{1,0\big\}$. The first fundamental form of $T_1S^2$ with respect to
    these coordinates is \cite{K-S}
    $$
    d{\tilde s}^2=du^2+dv^2+2\cos u\,dv\,d\omega +d\omega^2.
    $$
    The local parameterization of the submanifold $\xi(S^2)$, generated by the given field,
    is $ \omega =v$ and the induced metric on $\xi(S^2)$ is
    $$
    d{\tilde s}^2=du^2+2(1+ \cos u)\, dv^2=  du^2+4\cos^2u/2 \,dv^2.
    $$
    Thus, $\xi(S^2)$ is locally isometric to sphere $S^2$ of radius 2. Since
    $T_1S^2\stackrel{isom}{\approx}RP^3$ and there are no other totally geodesic submanifolds in
    $RP^3$ except $RP^2$, we see that $\xi(S^2)$ is a part of $RP^2$ .
    So the proof is complete.
    \end{proof}

    \begin{proposition}\label{Zero}
    Let $M^2$ be a Riemannian manifold of constant zero  curvature $K=0$. Then
    a totally geodesic unit vector field $\xi$ on $M^2$ is either parallel or
    moves along the family of parallel geodesics with  constant angle speed.
    Geometrically, $\xi(M^2)$ is either $E^2$ imbedded isometrically
    into $E^2\times S^1$ as a factor or a helical flat submanifold in $E^2\times S^1$.
    \end{proposition}

    \begin{proof}
    Suppose $\xi$ is totally geodesic unit vector field on $M^2$.
    Set $\Omega\equiv 0$ in Lemma \ref{Forms}. Then $\lambda\mu\equiv 0$. If
    $\lambda\equiv 0$ over some domain $D\subset M^2$, then $\xi $ is {\it parallel in this
    domain}.

    Suppose $\lambda\not\equiv 0$ in a domain $D\subset M^2$. Then
    $\mu\equiv 0$ on at least a subdomain $D'\subset D$. This means that the $e_0$-curves
    are geodesics in $D'$ and the field $\xi$ is parallel along them. Choose a family
    of $e_0$-curves and the orthogonal trajectories as a local coordinate net in $D'$.
    Then the first fundamental form of $M^2$ takes the form
    $\ds
    ds^2=du^2+f^2\,dv^2
    $
    and since $M^2$ is of zero curvature, $f$ satisfies the equation
    $$
    f_{uu}=0.
    $$
    A general solution of this equation is
    $\ds
    f(u,v)=A(v)u+B(v).
    $
    There are two possible cases:
    $$
    \begin{array}{rl}
        \mbox{(a)}  & A(v)\ne 0 \mbox{ in some subdomain } D''\subset D';\\
        \mbox{(b)}  & A(v)\equiv 0 \mbox{ over the whole domain } D'.
    \end{array}
    $$
    Case(a).
    The  function $f$ may be presented over $D''$ in the form
    $$
    f(u,v)=A(v)(u+\theta),
    $$
    where $\theta(v)=B(v)/A(v)$. After a $v$-parameter change, the
    metric in $D''$ takes the form
    $\ds
    ds^2=du^2+(u+\theta)^2\,dv^2.
    $
    Applying Proposition \ref{Semigeo} for $f=u+\theta$, we get
    $\ds
    \lambda=\frac{\omega'+1}{u+\theta}.
    $
    Setting $\Omega_{11}\equiv 0$, we obtain the identity
    $$
    \omega''(u+\theta)-(\omega'+1)\theta'\equiv 0.
    $$
    From this we get $\left\{\begin{array}{l} \omega''=0\\
    \omega'=-1\end{array}\right.$ or $\left\{\begin{array}{l} \omega''=0\\
    \theta'=0\end{array}\right.$. In the first case, $\lambda=0$ and the field $\xi$ is
    parallel again. In the second case
    $\ds
    \left\{
    \begin{array}{l}
        \theta=const, \\
        \omega=av+b
    \end{array}
    \right.
    $
    \ $ a,b=const$.

    Making a parameter change, we reduce the metric to the form
    $$
    ds^2=du^2+ u^2\,dv^2
    $$
    Applying Proposition \ref{Semigeo} with $f(u,v)=u$, we get
    $\ds
    \lambda = \frac{a+1}{u}.
    $
    The substitution into $\Omega_{01}$ gives the condition
    $$
       -\frac{a+1}{u^2+(a+1)^2}=0
    $$
    which is possible only if $a=-1$. But this means that again $\lambda=0$
    and hence $\xi$ is a parallel vector field.

    Case (b). After a $v$-parameter change, the metric takes the form
    $$
    ds^2=du^2+dv^2.
    $$
    Applying Proposition \ref{Semigeo} for $f\equiv 1$, we get
    $\ds
       \lambda=\omega'.
    $
    Setting $\Omega_{11}\equiv 0$, we find  $ \omega''\equiv 0$. This means that $\omega = av+b$ and $\xi $ is either
    parallel along the $u$-lines $(a=0)$ or moves along the $u$-lines helically with constant angle
    speed.

    Let $(u,v,\omega)$ be standard coordinates in $E^2\times S^1$. Then the first
    fundamental form of $E^2\times S^1$ is
    $$
    d{\tilde s}^2=du^2+dv^2+d\omega^2.
    $$
    If $a=0$, then with respect to these coordinates the local parameterization of
    $\xi(E^2)$  is $\omega=const$ and $\xi(E^2)$ is nothing else but $ E^2$ isometrically
    imbedded into $E^2\times S^1$. If $a\ne 0$, then the local parameterization of $\xi(E^2)$ is
    $\omega=av+b$ and the induced metric is
    $$
    d{\tilde s}^2= du^2+(1+a^2)\,dv^2
    $$
    which is flat. The imbedding is helical in the sense that this submanifold meets each
    flat element of the cylinder $p:E^2\times S^1\to S^1$ under constant angle
    $\varphi=\arccos \frac{ 1}{\sqrt{1+a^2}}$.
      So the proof is complete.
    \end{proof}

    \subsection{The curvature}

    The main goal of this section is to obtain an explicit formula for the Gaussian
    curvature of $\xi(M^2)$ and apply it to some specific cases. The first step is
    the following lemma.

    \begin{lemma}\label{Sec2}
     Let $\xi$ be a unit vector field on a 2-dimensional Riemannian manifold of Gaussian curvature $K$.
     In   terms of Lemma \ref{prop1},
     the sectional curvature $K_{T_1M}(\xi)$ of $T_1M$ along 2-planes tangent to $\xi(M)$ is given by
     $$
 K_{T_1M}(\xi)=\frac{K^2}{4}+\frac{K(1-K)}{1+\lambda^2}+(-1)^{s+1}\frac{\lambda}{1+\lambda^2}e_0(K).
     $$
     \end{lemma}
     \begin{proof}
     Let $\tilde \pi$ be a 2-plane tangent to $\xi(M)$. Then $\tilde X = e_0^h$ and
     $\tilde Y=\frac{1}{\sqrt{1+\lambda^2}}(e_1^h+\lambda\eta^v)$ form an orthonormal
     basis of $\tilde\pi$. So we may apply (\ref{Sec}) setting $X_1=e_0$, $X_2=0$,
     $Y_1=\frac{1}{\sqrt{1+\lambda^2}}e_1$, $ Y_2=\frac{\lambda}{\sqrt{1+\lambda^2}}\eta$.

     We get
     $$
     \begin{array}{l}\ds
     \big<R(X_1,Y_1)Y_1,X_1\big>=
     \frac{1}{1+\lambda^2}\big<R(e_0,e_1)e_1,e_0\big>=
     \frac{1}{1+\lambda^2}K, \\[2ex]\ds
     \|R(X_1,Y_1)\xi\|^2=\frac{1}{1+\lambda^2}\|R(e_0,e_1)\xi\|^2
     =\frac{1}{1+\lambda^2}K^2, \\[2ex]\ds
     \|R(\xi,Y_2)X_1\|^2=\frac{\lambda^2}{1+\lambda^2}\|R(\xi,\eta)e_0\|^2=
     \frac{\lambda^2}{1+\lambda^2}K^2, \\[2ex]\ds
    \big<(\nabla_{X_1}R)(\xi,Y_2)Y_1,X_1\big>=\frac{\lambda}{1+\lambda^2}
     \big<(\nabla_{e_0}R)(\xi,\eta)e_1,e_0\big>=-(-1)^{s}\frac{\lambda}{1+\lambda^2}e_0(K),
     \end{array}
     $$
    where $K$ is the Gaussian curvature of $M$. Applying directly (\ref{Sec}) we
    obtain
    $$
    \begin{array}{rl}
    K_{T_1M}(\xi)&\ds =\frac{1}{1+\lambda^2}\left(K-\frac34
    K^2+\frac{\lambda^2K^2}{4}+(-1)^{s+1}\lambda e_0(K)\right)\\[2ex]
    &\ds =\frac{1}{1+\lambda^2}\left(K(1-K)+\frac{(1+\lambda^2)K^2}{4}+(-1)^{s+1}\lambda
    e_0(K)\right)\\ [2ex]
    &\ds=\frac{K^2}{4}+\frac{K(1-K)}{1+\lambda^2}+(-1)^{s+1}\frac{\lambda}{1+\lambda^2} e_0(K).
    \end{array}
    $$
    \end{proof}

    Now we have the following.

    \begin{lemma}\label{main}
        Let $\xi$ be a unit vector field on a 2-dimensional Riemannian manifold $M$.
        In terms of Lemma \ref{prop1},
        the Gaussian curvature $K_\xi$ of the hypersurface $\xi(M)\in T_1M$ is given
        by
        $$
         \begin{array}{ll}
        \ds K_\xi=\frac{K^2}{4}+\frac{K(1-K)}{1+\lambda^2}&\ds +(-1)^{s+1}\frac{\lambda}{1+\lambda^2}
        e_0(K)+\\[2ex]
        &\ds\frac12
        \mu e_1\left(\frac{1}{1+\lambda^2}\right)-\left((-1)^{s+1}\frac{K}{2}+
        \frac{e_0(\lambda)}{1+\lambda^2} \right)^2,
        \end{array}
       $$
       where $K$ is the Gaussian curvature of $M$.
    \end{lemma}
    \begin{proof}
    In our case, one can  easily reduce the formula (\ref{tensor}) to the form
    $$
    K_\xi=K_{T_1M}(\xi)+\det \Omega.
    $$

    Applying Lemma \ref{Forms}, we see that
    $$
    \begin{array}{ll}
    \ds\det\Omega=&\ds-\mu\frac{\lambda}{\sqrt{1+\lambda^2}}
    e_1\left(\frac{\lambda}{\sqrt{1+\lambda^2}}\right)-
    \left((-1)^{s+1}\frac{K}{2}+\frac{e_0(\lambda)}{1+\lambda^2} \right)^2= \\[2ex]
    &\ds-\frac12\mu e_1\left(\frac{\lambda^2}{1+\lambda^2}\right)-
    \left((-1)^{s+1}\frac{K}{2}+\frac{e_0(\lambda)}{1+\lambda^2} \right)^2 =\\[2ex]
    &\ds\frac12 \mu e_1\left(\frac{1}{1+\lambda^2}\right)-
    \left((-1)^{s+1}\frac{K}{2}+\frac{e_0(\lambda)}{1+\lambda^2} \right)^2.
    \end{array}
    $$
    Combining this result with Lemma \ref{Sec2}, we get what was claimed.
    \end{proof}

    As an application of Lemma \ref{main} we  prove the following theorems.

    \begin{theorem}\label{Extr}
         Let $M^2$ be a 2-dimensional Riemannian manifold of  Gaussian curvature $K$.
        Suppose that $\xi$ is a unit geodesic vector field on $M^2$. Then
        the submanifold $\xi(M^2)\subset T_1M^2$ has non-positive
        extrinsic curvature.
    \end{theorem}
    \begin{proof}
         By definition, the extrinsic curvature of a submanifold is the difference
        between the sectional curvature of the submanifold and the sectional curvature of
        ambient space along the planes, tangent to the submanifold. In our case , this is
        $\det\Omega$.
        If $\xi $ is a geodesic vector field, then we may choose $e_0=\xi$ and then
        $\mu=k=0$. Therefore,  for the extrinsic curvature we get
        $$
        -\left((-1)^{s+1}\frac{K}{2}+\frac{e_0(\lambda)}{1+\lambda^2} \right)^2\leq0 \, .
        $$
    \end{proof}

    \begin{theorem}\label{Const}
         Let $M^2$ be a space of constant Gaussian curvature $K$.
        Suppose that $\xi$ is a unit geodesic vector field on $M^2$. Then
        $\xi(M^2)$ has constant Gaussian curvature in one of the following cases:
        \begin{itemize}
            \item[(a)] $K=-c^2<0$
            and $\xi$ is a normal vector field for the family of horocycles on
             the hyperbolic 2-plane $L^2$. In this case $K_\xi=-c^2$ and therefore
                $\xi(L^2)$ is locally isometric to $L^2$;
            \item[(b)] $K=0$ and $\xi$ is a parallel vector field on $M^2$. In this case
                $K_\xi=0$ and $\xi(M^2)$ is  also locally isometric to $M^2$;
            \item[(c)] $K=1$ and $\xi$ is any (local) geodesic vector field on the standard
                sphere $S^2$. In this case $K_\xi=0$.
        \end{itemize}
    \end{theorem}

    \begin{proof}

        Since $\xi$ is geodesic, we may set $e_0=\xi$, $e_1=\eta, \ s=1$. Taking into account
         (\ref{Fr1}) and (\ref{Fr2}), we see that $\lambda=-\kappa=-\sigma$.
         Lemma \ref{prop1} (b)  gives $-K=-e_0(\sigma)+\sigma^2$. So the result
        of Lemma \ref{main} takes the form
        $$
         \begin{array}{ll}
        \ds K_\xi=& \ds \frac{K^2}{4}+\frac{K(1-K)}{1+\sigma^2}-\left(\frac{K}{2}-
        \frac{e_0(\sigma)}{1+\sigma^2} \right)^2=\\[2ex]
         & \ds \frac{K^2}{4}+\frac{K(1-K)}{1+\sigma^2}-\left(\frac{K}{2}-
        \frac{K+\sigma^2}{1+\sigma^2} \right)^2=\\[2ex]
        & \ds  \frac{K(1-K)}{1+\sigma^2}+\frac{K(K+\sigma^2)}{1+\sigma^2}-
        \left(\frac{K+\sigma^2}{1+\sigma^2} \right)^2= \\[2ex]
        & \ds K-\left(\frac{K+\sigma^2}{1+\sigma^2} \right)^2.
        \end{array}
       $$

       Suppose that $K_\xi$  is constant. Then the following cases should be
       considered:

   (a) $\ds\sigma=const\ne0$. This means that the orthogonal trajectories
    of the field $\xi$ consist of curves of constant curvature.
    With respect to this natural coordinate system, the metric of $M^2$
    takes the form
    $\ds
    ds^2=du^2+f^2\, dv^2.
    $
    Set $\sigma=-c$. Then the function $f$ should satisfy  the equation
    $$
    \frac{f_u}{f}=c
    $$
    the general solution of which is $f(u,v)=A(v)e^{cu}$. After $v$-parameter
    change we obtain metric of the form
    $$
    ds^2=du^2+e^{2cu}\, dv^2.
    $$
    So, the manifold $M^2$ is locally isometric to the hyperbolic 2-plane $L^2$ of
    curvature $-c^2$
    and  the field $\xi$ is a geodesic field of (internal or external) normals to the family of
    horocycles.

    (b) $\ds\sigma=0$. Then evidently $\xi$ is a parallel vector field
        and therefore the manifold $M^2$ is locally Euclidean which implies $K_\xi$=0.

   (c) $\ds\sigma$ is not constant. Then  $K_\xi$ is constant if $K=1$ only.
    So, $M^2$ is contained in a standard sphere $S^2$ and the curvature of $\xi(S^2)$ does not depend
    on $\sigma$. Thus, the field $\xi$ is any (local) geodesic vector field.
    Evidently, $K_\xi=0$ for this case.

     \end{proof}

The case $(a)$ of the Theorem \ref{Const} has an interesting generalization of the
following kind.
\begin{theorem}\label{Foli}
    Let $L^2$ be a hyperbolic 2-plane of curvature $-c^2$. Then $T_1L^2$
    admits a hyperfoliation with leaves of constant intrinsic curvature $-c^2$ and of  constant
    extrinsic curvature $-\frac{c^2}{4}$. The leaves are generated by unit vector fields
    making a constant angle with a pencil of parallel geodesics on $L^2.$
\end{theorem}
    \begin{proof}
        Consider $L^2$ with  metric $\ds ds^2=du^2+e^{2cu}\,dv^2$
        and a family of vector fields
        $$
        \xi_\omega=\cos\omega X_1+\sin\omega X_2 \quad (\omega=const),
        $$
        where
        $\ds
        X_1=\big\{1,0\big\},\quad X_2=\big\{0,e^{-cu}\big\}
        $
        are the unit vector fields.

        Since $\ds\nabla_{X_1}\xi_\omega=0$, we may set
        $\ds
        e_0=X_1,\quad e_1=X_2
        $
        end therefore we have
        $\ds
        \sigma=-c, \quad \lambda=c .
        $
        Then, setting $K=-c^2$ and $\lambda=c$ in Lemma \ref{main}, we get
        $$
        K_\xi=-c^2.
        $$
        The extrinsic curvature of $\xi(L^2)$ is also constant since
        $$
        \det\Omega=-\frac14 c^2.
        $$

        Now fix a point $P_\infty$ at infinity  boundary of $L^2$ and draw a
        pencil of parallel geodesics from $P_\infty$ through each point of $L^2$.
        Define a family  of submanifolds $\xi_\omega(L^2)$ for this pencil. Evidently,
        through each point $(p,\zeta)\in T_1L^2$  there passes only one
        submanifold of this family. Thus, a family of submanifods $\xi_\omega$ form a
        hyperfoliation on $T_1L^2$ of constant intrinsic curvature $-c^2$ and constant
        extrinsic curvature $-\frac{c^2}{4}$.

        Geometrically, $\xi_\omega(L^2)$ is a family of
        coordinate hypersurfaces $\omega=const$ in $T_1L^2$. Indeed, let $(u,v,\omega)$
        form a natural local coordinate system on $T_1L^2$. Then the metric of
        $T_1L^2$ has the form
        $$
        ds^2=du^2+2e^{2u}dv^2+2\,dv d\,\omega +d\,\omega^2.
        $$
        With respect to these coordinates, the coordinate hypersurface $\omega=const$
        is nothing else but $\xi_\omega(L^2)$ and the induced metric is
        $$
        ds^2=du^2+2e^{2cu}dv^2.
        $$
         Evidently, its Gaussian curvature is constant and equal to $-c^2$.

        \end{proof}

\vspace{1cm}

\noindent
Department of Geometry,\\
Faculty of Mechanics and Mathematics,\\
Kharkiv National University,\\
Svobody Sq. 4,\\
 61077, Kharkiv,\\
Ukraine.\\
e-mail: yamp@univer.kharkov.ua

 \end{document}